\documentclass[11pt]{amsart}
\usepackage[active]{srcltx}
\usepackage{amsmath,amssymb}
\usepackage{latexsym}
\usepackage{color}
\usepackage{dsfont}
\usepackage[normalem]{ulem}
\usepackage{epsfig}
\usepackage{pgf}
\newtheorem{teo}{Theorem}[section]

\newtheorem{prop}[teo]{Proposition}

\begin{document}

\title{Blow-up for a fully fractional heat equation}

\author{R. Ferreira and A. de Pablo}

\address{Ra\'{u}l Ferreira
\hfill\break\indent  Departamento de An\'alisis Matem\'atico y Matem\'{a}tica Aplicada,
\hfill\break\indent Universidad Complutense de Madrid,
\hfill\break\indent
 28040 Madrid, Spain.
\hfill\break\indent  e-mail: {\tt raul$_-$ferreira@mat.ucm.es}}

\address{Arturo de Pablo
\hfill\break\indent  Departamento de Matem\'{a}ticas,
\hfill\break\indent Universidad Carlos III de Madrid,
\hfill\break\indent
 28911 Legan\'{e}s, Spain.
\hfill\break\indent  e-mail: {\tt arturop@math.uc3m.es}}

\

\begin{abstract}
We study the existence and behaviour of blowing-up solutions to the fully fractional heat equation
$$
\mathcal{M} u=u^p,\qquad x\in\mathbb{R}^N,\;0<t<T
$$
with $p>0$, where $\mathcal{M}$ is a nonlocal operator given by a space-time kernel $M(x,t)=c_{N,\sigma}t^{-\frac N2-1-\sigma}e^{-\frac{|x|^2}{4t}}\mathds{1}_{\{t>0\}}$, $0<\sigma<1$. This operator coincides with the fractional power  of the heat operator, $\mathcal{M}=(\partial_t-\Delta)^{\sigma}$ defined through semigroup theory. We characterize the global existence exponent $p_0=1$ and the Fujita exponent $p_*=1+\frac{2\sigma}{N+2(1-\sigma)}$, and study the rate at which the blowing-up solutions below $p_*$ tend to infinity, $\|u(\cdot,t)\|_\infty\sim (T-t)^{-\frac\sigma{p-1}}$.

\end{abstract}

\maketitle



\section{Introduction}

\label{sect-introduction} \setcounter{equation}{0}

We consider non-negative bounded solutions to the following nonlocal problem with memory
\begin{equation}\label{eq.principal}
\left\{
\begin{array}{ll}
\mathcal{M} u=u^p,\qquad & (x,t)\in \mathbb R^N\times (0,T),\\
u(x,t)=f(x,t),\qquad & (x,t)\in \mathbb R^N\times (-\infty,0],
\end{array}\right.
\end{equation}
where $N\ge1$, $p>0$, $f$ a given regular nonnegative function,
and where $0<T\le\infty$ is the maximal time of existence. The nonlocal operator $\mathcal{M}$ is defined by the  formula
\begin{equation}\label{master-def}
\mathcal{M} u(x,t)=\int_{-\infty}^t\int_{\mathbb{R}^N}\left(u(x,t)-u(z,s)\right)M(x-z,t-s)\,dzds,
\end{equation}
where the space-time kernel is given by
\begin{equation}\label{kernel-def}
\begin{array}{l}
M(x,t)=K_t(x)L(t),\\ [2mm]
 K_t(x)=(4\pi t)^{-\frac N2}e^{-\frac{|x|^2}{4t}},\quad L(t)=L_\sigma(t)=\dfrac{1}{|\Gamma(-\sigma)|}t^{-1-\sigma}.
 \end{array}
\end{equation}
Notice that the above integral is well defined provided $u$ is parabolic H\"older continuous $C_p^{2\sigma+\varepsilon}(\mathbb R^N\times (-\infty,T))$ for some $\varepsilon>0$, see \eqref{space-C-gamma}.
The equation in \eqref{eq.principal} is an example of a \emph{master equation} as defined in \cite{Metzler00}, it describes a continuous time random walk where the random jumps occur with
random time lag. The operator is nonlocal in space and time, since
the value of $\mathcal{M} u$ at a given point $(x,t)$  depends on the values of $u$ in the whole $\mathbb{R}^N$ and for all previous times $s<t$.

\

{\sc The operator $\mathcal{M}$}.
To look more closely to the operator we assume first $T=\infty$ so we can use Fourier transform $\mathcal{F}$ in both $x$ and $t$. We obtain
$$
\mathcal{F}(\mathcal{M} u)(\xi,\tau)=(i\tau+|\xi|^2)^{\sigma}\mathcal{F}(u)(\xi,\tau),
$$
see for instance \cite{Sampson68}. Since $m(\xi,\tau)=i\tau+|\xi|^2$ is precisely the symbol of the heat operator $\partial_t-\Delta$, our operator is, in Fourier sense, $\mathcal{M}=(\partial_t-\Delta)^{\sigma}$.
This is well defined for functions in the space
$$
\mathcal{H}^\sigma=\{u\in L^2(\mathbb{R}^{N+1})\,:\,m^{\frac\sigma2}\mathcal{F}(u)\in L^2(\mathbb{R}^{N+1})\}.
$$
We recall that we are mainly interested in solutions to equation \eqref{eq.principal} that are defined only up to a finite time, $t\in(-\infty,T)$, they are bounded in that interval and satisfy
$$
\lim_{t\to T^-}\|u(\cdot,t)\|_\infty=\infty.
$$
We say in that case that $u$ blows up at $t=T$.
We are thus not allowed to use Fourier transform in time, and our operator $\mathcal{M}$ \emph{is not} $(\partial_t-\Delta)^{\sigma}$, at least in Fourier sense, but when applied to globally defined functions.

On the other hand, the fractional power of the operator $\mathcal{L}=\partial_t-\Delta$ can be defined using the theory of semigroups by means of the formula
\begin{equation}\label{semigroups}
  \mathcal{L}^\sigma=\frac1{|\Gamma(-\sigma)|}\int_0^\infty\left(I-e^{-s\mathcal{L}}\right)\,\frac{ds}{s^{1+\sigma}},
\end{equation}
see \cite{Bochner}, where $v=e^{-s\mathcal{L}}u$ is the solution to the problem
\begin{equation}\label{dobleheat}
\left\{
\begin{array}{ll}
  \partial_sv+\mathcal{L}v=0,&s>0,\\ [2mm]
  v\big|_{s=0}=u.
\end{array}
\right.
\end{equation}
We then have
\begin{equation*}\label{dobleheat2}
\left\{
\begin{array}{ll}
  (\partial_s+\partial_t)v=\Delta v,&x\in\mathbb{R}^N,\;t,s>0,\\ [2mm]
  v(x,t,0)=u(x,t),&x\in\mathbb{R}^N,\; t>0,
\end{array}
\right.
\end{equation*}
and then $v$ is given, translated in time, by convolution with the Gauss kernel,
$$
v(x,t,s)=\int_{\mathbb{R}^N}K_{s}(z)u(x-z,t-s)\,dz.
$$
Introducing this expression into \eqref{semigroups}, using that $
\int_{\mathbb{R}^N}K_s(z)\,dz=1$, we obtain
$$
\begin{array}{rl}
\mathcal{L}^\sigma u(x,t)&\displaystyle=\frac1{|\Gamma(-\sigma)|}\int_0^\infty\left(u(x,t)-v(x,t,s)\right)\,
\frac{ds}{s^{1+\sigma}}, \\ [3mm]
&\displaystyle=\frac1{|\Gamma(-\sigma)|}\int_0^\infty\int_{\mathbb{R}^N}\left(u(x,t)-u(x-z,t-s)\right)K_s(z)\,dz\,
\frac{ds}{s^{1+\sigma}},
\end{array}
$$
which coincides with \eqref{master-def}. Since this calculation does not depend on $u$ being defined for all positive times, and only a minimum of regularity is needed in order for the integral to be well defined, it can be concluded that in the sense of semigroups our operator is indeed $\mathcal{M}=(\partial_t-\Delta)^{\sigma}$.

The operator $(\partial_t-\Delta)^{\sigma}$ is studied in detail in \cite{StingaTorrea17}, see also \cite{Athana-eta18,NystromSande16,Taliaferro21}. We recall the main immediate properties, valid also for the operator $\mathcal{M}$: it is linear, invariant under translations in space and time, and homogeneous of order $2\sigma$ under scaling $x\to\lambda x$, $t\to\lambda^2 t$; a (strong) maximum principle holds. More properties are obtained in \cite{Athana-eta18,StingaTorrea17} by means of the extension technique, see the next section.

When applied to functions depending only on one variable, space or time, the operator $\mathcal{M}$ simplifies. In fact we have, for $u=u(t)$,
$$
\begin{array}{rl}
\mathcal{M} u(t)&\displaystyle=\int_{-\infty}^t\left(u(t)-u(s)\right)L(t-s)\int_{\mathbb{R}^N}K_{t-s}(x-z)\,dzds
\\ [3mm]
&\displaystyle=\int_{-\infty}^t\left(u(t)-u(s)\right)L(t-s)\,ds.
\end{array}$$
This is the Marchaud derivative, see for instance \cite{Samko78},
\begin{equation}\label{Marchaud}
\partial_t^\sigma u(t)=\frac{1}{|\Gamma(-\sigma)|}\int_{-\infty}^t\frac{u(t)-u(s)}{(t-s)^{1+\sigma}}\,ds.
\end{equation}
Also, for $u=u(x)$, since
$$
\int_0^{\infty} t^{-\frac N2-1-\sigma}e^{-\frac{|x|^2}{4t}}\,dt=4^{\frac N2-1+\sigma}\Gamma(\frac N2+\sigma)
|x|^{-N-2\sigma},
$$
we obtain
$$
\begin{array}{rl}
\mathcal{M} u(x)&\displaystyle=\int_{\mathbb{R}^N}\left(u(x)-u(z)\right)\int_{-\infty}^t M(x-z,t-s)\,dzds \\ [3mm]
&\displaystyle=c_{N,\sigma}\int_{\mathbb{R}^N}\frac{u(x)-u(z)}{|x-z|^{N+2\sigma}}\,dz=(-\Delta)^\sigma u(x),
\end{array}
$$
the well-known fractional Laplacian, \cite{Landkof72}, with the exact constant, see for example \cite{dPQRV}. In other words,
\begin{equation}\label{separate}
(\partial_t-\Delta)^\sigma u(t)=\partial_t^\sigma u(t),\qquad (\partial_t-\Delta)^\sigma u(x)=(-\Delta)^\sigma u(x).
\end{equation}

On the very other hand, $\mathcal{M}$  has a fundamental solution
\begin{equation}\label{fundamental}
G(x,t)=K_t(x)L_{-\sigma}(t)\mathds{1}_{\{t>0\}}=A_{N,\sigma}t^{-\frac N2-1+\sigma}e^{-\frac{|x|^2}{4t}}\mathds{1}_{\{t>0\}},
\end{equation}
satisfying
$$
\mathcal{M}G(x,t)=\delta(x,t).
$$
We just have to apply formula \eqref{master-def} (or formally calculate the  Fourier transform of $G$, which is $m^{-\sigma}(\xi,\tau)$). In our particular case we have a representation formula for the solution to the memory value problem
\begin{equation}\label{eq.general}
\left\{
\begin{array}{ll}
\mathcal{M} u=h,\qquad & \text{for } t>0,\\
u=f,\qquad & \text{for } t\le0,
\end{array}\right.
\end{equation}
which if $h$ and $f$ are regular enough (see the precise conditions in the Appendix) is
\begin{equation}\label{representation}
u(x,t)=\int_{-\infty}^t\int_{\mathbb{R}^N}F(z,s)G(x-z,t-s)\,dzds,
\end{equation}
with
\begin{equation}\label{representation-F}
F=\begin{cases}
  h&\text{if } t>0, \\
  \mathcal{M}f&\text{if } t\le0.
\end{cases}
\end{equation}
Thus the solution to our problem \eqref{eq.principal} can be written in the form
\begin{equation}\label{repres-I12}
\begin{array}{rl}
u(x,t)&\displaystyle=\int_{-\infty}^0 \int_{\mathbb{R}^N}\mathcal{M}f(z,s)G(x-z,t-s)\,dzds,\\  [3mm]
&\displaystyle +\int_{0}^t\int_{\mathbb{R}^N}u^p(z,s)G(x-z,t-s)\,dzds.
\end{array}
\end{equation}
We say that $u$ is a mild solution to problem \eqref{eq.principal} if \eqref{repres-I12} holds in $\mathbb{R}^N\times(-\infty,T)$. Observe that the second term only requires $u$ bounded up to $t$ in order to be well defined. As to the first term, since the Green function is not integrable in time at $-\infty$, we impose the following hypothesis on the memory data $f$,
\begin{equation}\label{H}
\tag{H}
|\mathcal{M}f(z,s)|\le c|s|^{-\gamma},\quad\text{for } s\le-1\text{ uniformly in } x\in\mathbb{R}^N,
\end{equation}
for some $\gamma>\sigma$. This condition is satisfied for instance if $f$ is $C^1$ in time and
both $f$ and $\partial_tf$ decay as $|s|^{-\gamma}$.

We will see in the Appendix that mild solutions are classical if $f$ is regular enough.

\

{\sc Previous results on blow-up}

In the case $\sigma=1$ problem \eqref{eq.principal} reduces to the standard semilinear equation
\begin{equation}\label{local-eq}
\partial_tu-\Delta u=u^p,
\end{equation}
a well-known equation studied since the sixties of the previous century in relation to blow-up with the works of Kaplan and Fujita, \cite{Fujita,Kaplan}. For that equation there exist two critical exponents, the global existence exponent $p_0=1$ and the \emph{Fujita} exponent $p_*=1+\frac2N$, such that: $(i)$ all solutions are global in time if $p\le p_0$, $(ii)$ all solutions blow-up in finite time if $p_0<p\le p_*$, and $(iii)$ there exist both, global in time solutions and blow-up solutions if $p>p_*$. From those results a lot of research has been performed for different local operators, like the $p$--Laplacian or the Porous Medium. See for instance the review books \cite{QuittnerSouplet07,GKMS}.

In the fractional derivatives framework we quote the works \cite{KLT05,NagasawaSirao69,ZS15} for the space fractional, time fractional or even double fractional problem
\begin{equation}\label{doublefrac}
\left\{
\begin{array}{ll}
D^\alpha_tu+(-\Delta )^{\beta}u=u^p&\quad x\in\mathbb{R}^N,\;t>0, \\ [2mm]
u(x,0)=f(x)&\quad x\in\mathbb{R}^N,
\end{array}\right.
\end{equation}
with $0<\alpha\le1$ and $0<\beta\le1$. Here $D^\alpha_t$ is the Caputo derivative (Marchaud derivative when considering constant data $u(x,t)=u(x,0)$ for every $t<0$); as far as we know there are no works on blow-up dealing with Marchaud derivative. The papers \cite{NagasawaSirao69} and \cite{ZS15} consider, respectively, the cases $\alpha=1$ and $\beta=1$, showing that the (Caputo) Fujita exponent is $p_C=1+\frac{2\beta}N$. In the case $0<\sigma<\beta=1$, it is also proved that for the critical exponent $p=1+\frac2N$ there exist global solutions, contrary to what happens for $\alpha=1$.

The inner case $0<\alpha,\beta<1$ has been studied in \cite{KLT05}, where the authors prove the estimate
\begin{equation}\label{fujita-KLT}
p_C\ge1+\frac{2\alpha\beta}{N\alpha+2\beta(1-\alpha)},
\end{equation}
which in view of the previous result for $\beta=1$ seems not to be sharp.

On the other hand, Taliaferro has recently proved in \cite{Taliaferro21} that the inequality problem
\begin{equation}\label{Taliaferro}
\left\{
\begin{array}{ll}
(\partial_t-\Delta)^\sigma u\ge u^p&\quad x\in\mathbb{R}^N,\;t\in \mathbb{R}, \\ [2mm]
u(x,t)=0&\quad x\in\mathbb{R}^N,\;t<0,
\end{array}\right.
\end{equation}
possesses no nontrivial solution $u\in L_{loc}^r((0,\infty);L^r(\mathbb {R}^N))$, with $1\le r<\frac{p(N+2)}{2\sigma}$, whenever
\begin{equation}\label{fujita-taliaferro}
0<p\le p_T=1+\dfrac{2\sigma}{N+2(1-\sigma)},\quad p\neq1.
\end{equation}
Though no conclusion can be deduced from \cite{Taliaferro21} in relation to our problem \eqref{eq.principal}, observe  that the exponent in~\eqref{fujita-taliaferro} coincides with the exponent estimate \eqref{fujita-KLT} obtained in \cite{KLT05} when putting $\alpha=\beta=\sigma$, and thus it could serve as an interesting conjecture.

\

{\sc Fujita exponent.}
Then, what is the expected Fujita exponent for problem \eqref{eq.principal}? There are two homogeneity arguments usually applied to similar reaction-diffusion problems to guess the answer: $(i)$ the equilibrium of the diffusion exponent decay with the ODE blow-up exponent, see \cite{DengLevine00}; and $(ii)$ the invariance of some $L^q$ norm, $q>1$, under rescaling, see \cite{CazenaveDW08}.

In the semilinear heat equation \eqref{local-eq}, the diffusion decay rate is $\rho_d=\frac N2$ (the decay of the Gauss kernel), while the solutions of the corresponding ODE $\partial_tz=z^p$ are $z(t)=c(T-t)^{-\frac1{p-1}}$, so $\rho_r=\frac1{p-1}$. Thus $\rho_d=\rho_r$ gives $p*=1+\frac2N$.
The same can be applied for the double fractional problem \eqref{doublefrac}, but only if $\alpha=1$ or $2\beta>N$, since otherwise the fundamental solution is not bounded at $x=0$, see \cite{Kemppainen-etal16}. If $\alpha=1$ the fractional Gauss kernel decays like $t^{-\frac{N}{2\beta}}$; if $0<\alpha<1=N<2\beta$, the fundamental solution decays like $t^{-\frac\alpha{2\beta}}$; the solution to $D_t^\alpha z=z^p$ satisfies $z(t)\sim(T-t)^{-\frac\alpha{p-1}}$. We then have in both cases $\rho_d=\frac{N\alpha}{2\beta}$, $\rho_r=\frac\alpha{p-1}$, producing $p_C=1+\frac{2\beta}N$. This again suggests that \eqref{fujita-KLT} is not sharp.
Now for problem \eqref{eq.principal}, and in view of \eqref{separate} (the solution to $\partial_t^\sigma z=z^p$ is $z(t)=c(T-t)^{-\frac\sigma{p-1}}$), and using the decay of the Green function \eqref{fundamental}, the exponents are
$$
\rho_d=\frac N2+1-\sigma,\quad \rho_r=\frac\sigma{p-1}\quad\rightsquigarrow\quad p_*=1+\dfrac{2\sigma}{N+2(1-\sigma)},
$$
exactly that in \eqref{fujita-taliaferro}.

On the other hand, the author in \cite{Weissler} proved for problem \eqref{local-eq} that if $\|u_0\|_{q_c}$ is small, where $q_c=\frac{N(p-1)}2>1$, then the solution is global in time. This exponent $q_c$ is the only one for which the above norm remains invariant under the natural scaling $u_\lambda(x,t)=\lambda^{\frac2{p-1}}u(\lambda x,\lambda^2 t)$ of the equation. The condition $q_c>1$ implies $p>1+\frac2N=p_*$. Applied this numerical argument to the equation \eqref{doublefrac} gives $q_c=\frac{N(p-1)}{2\beta}$ and thus the conjecture $p_C=1+\frac{2\beta}N$, proved to be correct if $\alpha=1$ or $\beta=1$. Since the homogeneity of the equations in \eqref{eq.principal} and in \eqref{doublefrac} is the same, this same calculus will produce for problem \eqref{eq.principal} an exponent $p_*=1+\frac{2\sigma}N$, different from the previously guessed. We can see that this is not the correct argument for the \emph{memory} problem \eqref{eq.principal}. Instead, the $L^{q_c}$ norm must be calculated both in space and time, and also, suggested by the representation formula \eqref{repres-I12}, it is the data $\mathcal{M}f$ that must be taken into consideration, not just $f$. We have the rescaling $u_\lambda(x,t)=\lambda^{\frac{2\sigma }{p-1}}u(\lambda x,\lambda^2 t)$, while
$$
\int_{-\infty}^0 \int_{\mathbb{R}^N}\left|\mathcal{M}f_\lambda(x,t)\right|^q\,dxdt=
\lambda^{\frac{2\sigma qp}{p-1}-N-2}
\int_{-\infty}^0 \int_{\mathbb{R}^N}\left|\mathcal{M}f(x,t)\right|^q\,dxdt.
$$
This gives
$$
q_c=\frac{(N+2)(p-1)}{2p\sigma}>1\quad\rightsquigarrow\quad p>1+\frac{2\sigma}{N+2(1-\sigma)}.
$$
We show in this paper that this is indeed the Fujita exponent. The global existence exponent is easy to derive, and since the operator is linear it has to do with the reaction term being sublinear or superlinear.

\begin{teo}\label{teo-BUP}
\begin{itemize}
  \item If $p\le1$  then the solution to problem \eqref{eq.principal} is global in time.
  \item If $1<p\le 1+\frac{2\sigma}{N+2(1-\sigma)}$ then the solution $u$  to problem \eqref{eq.principal} blows up in finite time.
      \item If $p> 1+\frac{2\sigma}{N+2(1-\sigma)}$ then there exist global solutions to problem \eqref{eq.principal}, and also blowing-up solutions, depending on the memory data $f$.
\end{itemize}
\end{teo}

To prove those results we follow the techniques used for the analogous local problem \eqref{local-eq} with the modifications needed to deal with nonlocal operators. Also we consider for some proofs an equivalent  local problem obtained when adding an extra variable, see \cite{NystromSande16,StingaTorrea17} and the next section, in the spirit of the well known Caffarelli-Silvestre extension \cite{CaffarelliSilvestre07} for the fractional Laplacian.

\

{\sc Blow-up rates.}
We are also interested in determining the velocity at which the nonglobal solutions blow up, the so-called \emph{blow-up rate}. For the local equation \eqref{local-eq} it is proved, under some restrictions, that the generic behaviour is
\begin{equation}\label{local-rates}
\|u(\cdot,t)\|_\infty\sim (T-t)^{-\frac1{p-1}}.
\end{equation}

We show here a generic blow-up rate with the exponent given by the time evolution equation $\partial_t^\sigma z(t)=z^p(t)$, as commented upon before. We thus prove that the blow-up rate is $\frac\sigma{p-1}$. We apply a very simple rescaling argument that works below the Fujita exponent.

\begin{teo}\label{teo-rate}
If $1<p\le1+\frac{2\sigma}{N+2(1-\sigma)}$ then the solution to problem \eqref{eq.principal} satisfies
\begin{equation}\label{rate}
\|u(\cdot,t)\|_\infty\sim (T-t)^{-\frac\sigma{p-1}}
\end{equation}
for $t\sim T$.\end{teo}

\

{\sc Organization of the paper.} After presenting the mentioned extension technique in Section~\ref{sec-extended}, we prove Theorem~\ref{teo-BUP} in Section~\ref{sec-bup}. Finally Theorem~\ref{teo-rate} is proved in Section~\ref{sec-rate}. The existence of solution and some properties, Theorem~\ref{teo-existence}, is
left to an Appendix.

Throughout the paper the letters $c,c_i,\,$ denote different constants not depending on the relevant quantities.

\

\section{Extended problem}\label{sec-extended}
Following the ideas of Stein \cite{Stein70} and Caffarelli and Silvestre \cite{CaffarelliSilvestre07},
given a function $u:\mathbb{R}^N\times(-\infty,T)\to\mathbb{R}$, in \cite{NystromSande16,StingaTorrea17} it is defined the $\sigma$--caloric extension $U$ of $u$,
$$
U=E(u):\Lambda\times(-\infty,T)\to\mathbb{R},
$$
to the upper half-space $\Lambda=\mathbb{R}^N\times\mathbb{R}_+$ by means of the formula,
$$
U(X,t)\equiv U(x,y,t)=\int_{-\infty}^t\int_{\mathbb{R}^N}u(z,s)P_y(x-z,t-s)\,dzds,
$$
where $P_y$ is the  kernel
$$
P_y(x,t)=d_{N,\sigma}y^{2\sigma}t^{-\frac N2-1-\sigma}e^{-\frac{|x|^2+y^2}{4t}},
$$
with coefficient $d_{N,\sigma}^{-1}=(4\pi)^{\frac N2}2^{2\sigma}\Gamma(\sigma)$. Since $P_y$ is positive with
$$\int_0^{\infty}\int_{\mathbb{R}^N}P_y(x,t)\,dxdt=1,
$$
the extension is order preserving with $E(1)=1$.

The importance of this extension function relies on the following two properties proved in \cite{NystromSande16,StingaTorrea17}:

$(i)$ We can obtain the fractional operator $\mathcal{M}$ of \eqref{master-def} using the normal derivative of the extended function on $\partial\Lambda=\mathbb{R}^N$ through the formula
\begin{equation}\label{conormal}
\mathcal{M} u(x,t)=-\kappa_{\sigma}\lim_{y\to0^+}y^\gamma \partial_yU(x,y,t),
\end{equation}
where $\kappa_\sigma=\frac{d_{N,\sigma}A_{N,\sigma}}{2\sigma}=\frac{\Gamma(\sigma)}{2^{\gamma}\Gamma(1-\sigma)}$, $\gamma=1-2\sigma$. \newline
In fact this would hold with any function $R(x,y,t)$ as kernel provided $y^{-2\sigma}R\big|_{y=0}=cM$ and $\iint R=1$.

$(ii)$ More relevant, for that particular choice of $P$ the function $U$ satisfies the extended (local) heat problem in $\Lambda$
\begin{equation}\label{ext}
\left\{
\begin{array}{ll}
y^\gamma \partial_tU-\text{div}(y^\gamma \nabla U)=0,\qquad&X\in\Lambda,\;0<t<T,\\ [2mm]
U(x,0,t)=u(x,t),&x\in\mathbb{R}^N,\;-\infty<t<T.
\end{array}
\right.
\end{equation}
Here the operators $\text{div}$ and $\nabla $ refer to the variable $X=(x,y)$. Thus $P$ is in fact the Poisson kernel for problem~\eqref{ext}. Though the proof, for instance in \cite{StingaTorrea17}, uses Fourier transform and holds in the sense of $L^2(\mathbb{R}^{N+1})$, it makes sense and can be justified whenever the integrals involved are well defined.

With these two properties we can write problem \eqref{eq.principal} in the equivalent form
\begin{equation}\label{eq.principal2}
\left\{
\begin{array}{ll}
y^\gamma \partial_tU-\text{div}(y^\gamma \nabla U)=0,&X\in\Lambda,\;0<t<T,\\ [2mm]
\displaystyle-\kappa_{\sigma}\lim_{y\to0^+}y^\gamma \partial_yU(x,y,t)= U^p(x,0,t),&x\in\mathbb{R}^N,\;0<t<T,\\ [2mm]
U(X,0)=g(X),&X\in\Lambda,
\end{array}
\right.
\end{equation}
with $g=E(f)\big|_{t=0}$, that is
$$
g(x,y)=\int_{-\infty}^0\int_{\mathbb{R}^N}f(z,s)P_y(x-z,-s)\,dzds.
$$
Clearly $\|g\|_\infty\le \|f\|_\infty$.

We then recover our function $u$ by the identity $u(x,t)=U(x,0,t)$. In the sequel we use both problems, \eqref{eq.principal} and \eqref{eq.principal2}, at our best convenience.

Observe that the case $\gamma=0$ (i.e. $\sigma=1/2$), reduces problem \eqref{eq.principal2} to the heat equation in a half-space with Neumann reaction condition. This is a well studied problem, see for instance \cite{DengFilaLevine94}, where it is shown that the Fujita exponent is $p_*=1+\frac1{N+1}$.

The equation in \eqref{eq.principal2} is  a degenerate equation with a Muckenhoupt $A_2-$weight given by $A(x, y) = y^\gamma$, see for instance \cite{ChiarenzaSerapioni}. This allows to obtain
existence, uniqueness, comparison and regularity, see \cite{Athana-eta18,StingaTorrea17} for the corresponding linear problem. We then apply this result in order to consider the nonlinear problem \eqref{eq.principal}. For that purpose we need the data $f$ being H\"older continuous of order $\gamma>\sigma$ in time and $2\gamma$ in space, that we denote $f\in C^\gamma_p$, see the Appendix; we also impose the decreasing condition \eqref{H} in time.

\

\begin{teo}\label{teo-existence}
For every $f\in C_p^\gamma(\mathbb{R}^N\times(-\infty,0))$, satisfying hipothesis \eqref{H}  there exists a unique  solution $u\in C_p^{\gamma+2\sigma}(\mathbb{R}^N\times(0,T))$ to problem \eqref{eq.principal} defined in a maximal time of existence $T\le\infty$. Moreover a comparison principle holds: if $u_1,u_2$ are the solutions to problem \eqref{eq.principal} with data  $f_1,f_2$, respectively, then $f_1\le f_2$ implies $u_1\le u_2$. If $f\in C^\infty$ then $u\in C^\infty$.
\end{teo}

We also recall that the extension of the fundamental solution is
$$
\widetilde{G}(X,t)=E(G)(x,y,t)=e^{-\frac{y^2}{4t}}G(x,t).
$$
This is not a fundamental solution for problem \eqref{eq.principal2}, since this problem is not invariant under translations in the vertical variable. Nevertheless, the extended solution $U$ is given in terms of $\widetilde{G}$ by
\begin{equation}\label{representation-ext}
U(x,y,t)=\int_{-\infty}^t\int_{\mathbb{R}^N}F(z,s)\widetilde{G}(x-z,y,t-s)\,dzds ,
\end{equation}
where $F$ is given in \eqref{representation-F}. See \cite{StingaTorrea17}.

We finally observe that there is a nonincreasing quantity defined for the extended problem \eqref{eq.principal2},
\begin{equation}\label{energy}
I_U(t)=\frac12\int_\Lambda |\nabla U(X,t)|^2\,d\mu(X)-\frac{1}{(p+1)\kappa_\sigma}\int_{\mathbb{R}^N}U^{p+1}(x,0,t)\,dx,
\end{equation}
where $d\mu(X)=y^\gamma\,dxdy$.
This is well defined whenever $U(\cdot,t)\in H^1(\Lambda;d\mu)$. In fact if $U$ solves \eqref{eq.principal2} then
$$
I'_U(t)=-\int_\Lambda |\partial_tU(X,t)|^2\,d\mu(X)\le0.
$$
It will play a role in the characterization of global solutions in the range above the Fujita exponent.

\

\section{Blow-up}\label{sec-bup}
We prove in this section Theorem \ref{teo-BUP}.
We first show that the global exponent is $p_0=1$, that is, for $p\le1$ the solutions are global and bounded locally in time.

\begin{teo}\label{global}
If $p\le1$  then the solution to problem \eqref{eq.principal} is global in time.
\end{teo}

\begin{proof} Let first $p<1$ and assume that there exists a time $-t_0<0$ such that $f(x,t)\equiv0$ for $t\le -t_0$. We compare with the global solution
$$
\overline u(x,t)=c_*(t+t_1)_+^{\nu},\qquad \nu=\frac\sigma{1-p},\quad c_*^{p-1}=\frac{\Gamma(1+\nu)}{\Gamma(1+\nu-\sigma)}.
$$
Just observe that $\overline u$ solves the equation $\partial_t^\sigma z=z^p$, where $\partial_t^\sigma$ is the Marchaud derivative, see \eqref{Marchaud}. We only need to compare the data, and for that purpose we choose $t_1> t_0$ such that
$$
c_*(t_1-t_0)^\nu\ge \|f\|_\infty.
$$
In the case $p=1$ the function $\overline u(t)=e^{t+t_1}$ does the job provided $e^{t_1-t_0}\ge \|f\|_\infty$.

To consider the general case without the restriction of $f$ having compact support in time, we compare in the extended problem with a supersolution of the form
$$
\overline V(x,y,t)=e^{\alpha t}\varphi(y).
$$
The profile $\varphi$ must satisfy:
$$
\begin{array}{ll}
\displaystyle \alpha y^\gamma \varphi\ge (y^\gamma \varphi')',&\quad y>0,\\ [2mm]
\displaystyle-\kappa_{\sigma}e^{\alpha t}\lim_{y\to0^+}y^\gamma \varphi'(y)\ge e^{p\alpha t}\varphi^p(0),&\quad t>0,\\ [2mm]
\displaystyle e^{\alpha t}\varphi(0)\ge \|f\|_\infty,&\quad t>0.
\end{array}
$$
By a simple rescaling we can get rid of the coefficient $\kappa_\sigma$. Then
these three conditions can be fulfilled, after straightforward computations, with the function
$$
\varphi(y)=A(2-y^{2\sigma}e^{-y^2}),\qquad A \text{ and } \alpha \text{ large}.
$$

\end{proof}

We now construct blow-up solutions above $p_0$.

\begin{teo}\label{teo-kaplan}
If $p>1$ then the solution to problem \eqref{eq.principal} blows up in finite time if the datum $f$ is large enough.
\end{teo}

\begin{proof}
Here and in the next proof we follow the standard argument of \cite{Fujita} and \cite{Kaplan} for the local case using an appropriate test function, this time applied to the extended problem~\eqref{eq.principal2}. Let $V$ be a solution to problem \eqref{eq.principal2} with some initial value
$g$, a nonnegative given function. Assume also $\partial_yg\le0$ in $\Lambda$. This gives  $\partial_y U\le0$ in $\Lambda$ for every $t>0$, just by looking at the problem satisfied by $\partial_y U$ and a comparison argument.

Consider the function
$$
\varphi(X)=\varrho  e^{-|X|^2},\quad \varrho=\dfrac2{\pi^{\frac N2}\Gamma(1-\sigma)}.
$$
It is easy to check that $\varphi$ has (weighted) integral 1.
$$
\displaystyle\int_\Lambda\varphi(X) \,d\mu(X)=1.
$$
Also, applying the  diffusion operator in \eqref{eq.principal2}, we have
\begin{equation}\label{gauss}
\text{div}(y^\gamma \nabla\varphi)\ge -(2N+4(1-\sigma))y^\gamma\varphi.
\end{equation}
Put now $\Lambda_\varepsilon=\mathbb{R}^N\times(\varepsilon,\infty)$, and define
$$
J_\varepsilon(t)=\displaystyle\int_{\Lambda_\varepsilon} V(X,t)\varphi(X)\,d\mu(X).
$$
Using the equation and integrating by parts we get
$$
\begin{array}{rl}
J_\varepsilon'(t)&\displaystyle=\int_{\Lambda_\varepsilon}\partial_tV\varphi \,d\mu(X)=\int_{\Lambda_\varepsilon} \text{div}(y^\gamma \nabla V)\varphi\,dxdy \\ [3mm]
&\displaystyle=\int_{\mathbb{R}^N}V(x,\varepsilon,t)2\varrho \varepsilon^{2-2\sigma}e^{-(|x|^2+\varepsilon^2)}\,dx\\ [3mm]
&\displaystyle\quad-\rho  \int_{\mathbb{R}^N}e^{-(|x|^2+\varepsilon^2)}\varepsilon^{1-2\sigma}\partial_y V(x,\varepsilon,t)\,dx\\ [3mm]
&\displaystyle\quad+\int_{\Lambda_\varepsilon} V\text{div}(y^\gamma \nabla\varphi)\,dxdy
\end{array}
$$
Letting $\varepsilon\to0^+$, and using \eqref{gauss} and \eqref{conormal} we get, for $J=J_0$,
$$
J'(t)\ge c_1 \int_{\mathbb{R}^N}e^{-|x|^2}V^p(x,0,t)\,dx-c_2J(t).
$$
On the other hand, by Jensen inequality,
$$
\begin{array}{rl}
J^p(t)&\displaystyle\le\int_{\Lambda} V^p(X,t)\varphi(X)\,d\mu(X)\\ [3mm]
&\displaystyle=\int_0^\infty\int_{\mathbb{R}^N}\rho y^\gamma e^{-(|x|^2+y^2)}V^p(x,y,t)\,dxdy \\ [3mm]
&\displaystyle\le \int_{\mathbb{R}^N}e^{-|x|^2}V^p(x,0,t)\,dx\int_0^\infty\varrho y^\gamma e^{-y^2}\,dy\\ [3mm]
&\displaystyle=c\int_{\mathbb{R}^N}e^{-|x|^2}V^p(x,0,t)\,dx.
\end{array}
$$
Therefore $J$ satisfies the differential inequality
\begin{equation}\label{ODE-J}
J'(t)\ge c_1 J^p(t)-c_2J(t),
\end{equation}
which implies blow-up of $J(t)$ for every $p>1$ provided the initial value $J(0)=\int_{\Lambda} g\varphi \,d\mu$ is large. Now, since
$$
\|u(\cdot,t)\|_\infty\ge\|U(\cdot,t)\|_\infty\ge J(t),
$$
by comparison we get blow-up for problem \eqref{eq.principal} for initial data $f$  satisfying $E(f)\ge g$ and $\partial_y g\le0$.
\end{proof}

\begin{teo}\label{teo-fujita}
If $1<p<1+\frac{2\sigma}{N+2(1-\sigma)}$ then for any given data $f$ the solution to problem \eqref{eq.principal} blows up in finite time.
\end{teo}

\begin{proof}
  We  modify the previous proof by introducing the rescaled test function
$$
\varphi_k(X)=\varrho k^\eta  e^{-k|X|^2},
$$
$\eta=\frac{N+2-2\sigma}2$. Clearly again
$\int_\Lambda\varphi_k \,d\mu=1$, while \eqref{gauss} becomes here
\begin{equation}\label{gauss2}
\text{div}(y^\gamma \nabla\varphi_k)\ge -cky^\gamma\varphi_k.
\end{equation}
Introducing the function $J(t)=\int_\Lambda V\varphi_k\,d\mu(X)$, we obtain exactly as before the differential inequality
\begin{equation}\label{ODE-Jk}
J'(t)\ge c_1 k^{1-\sigma}J^p(t)-c_2kJ(t).
\end{equation}
This produces blow-up whenever $J(0)> ck^{\frac\sigma{p-1}}$, or which is the same,
$$
\int_\Lambda  g(X)\,d\mu(X)> ck^{\frac\sigma{p-1}-\eta}.
$$
This holds for every initial value $g$ taking $k$ small if $p<1+\frac{\sigma}{\eta}=p_*$. We conclude with a comparison argument as before blow-up for problem \eqref{eq.principal}.
\end{proof}

To deal with the critical case $p=p_*$ we need some previous estimate using the representation formula \eqref{repres-I12}.

\begin{prop} Any solution to problem \eqref{eq.principal}, as a supersolution to the homogeneous equation, satisfies
\begin{equation}\label{eq.homogeneous2}
u(x,t)\ge ct^{\sigma-1}K_{t}(x)
\end{equation}
for every $x\in\mathbb{R}^N$, $t>t_0$ large and some $c>0$.
\end{prop}

\begin{proof}
Using the representation formula \eqref{representation} we get
$$
u(x,t)\ge\int_{-\infty}^0 \int_{\mathbb{R}^N}\mathcal{M}f(z,s)G(x-z,t-s)\,dzds.
$$
By the strong maximum principle proved in \cite[Lemma~1.5]{StingaTorrea17} we get $u>0$. Now, though $\mathcal{M}(f)$ has not a sign, the property that the integral is strictly positive, together with the fact that $G$ is also positive, implies that there exists an interval of negative times $[-b,-a]$ such that (by translation in space)
$$
\mathcal{M}f(z,s)\ge\varepsilon>0\qquad \text{for a.e. } |z|<\delta,\quad s\in[-b,-a],
$$
and
$$
\begin{array}{l}
\displaystyle\int_{-b}^{-a} \int_{|z|<\delta}\mathcal{M}f(z,s)G(x-z,t-s)\,dzds \\ [3mm]
\hspace{3cm}\displaystyle\le \frac12 \int_{-\infty}^0 \int_{\mathbb{R}^N}\mathcal{M}f(z,s)G(x-z,t-s)\,dzds.
\end{array}
$$
Thus
$$
\begin{array}{rl}
u(x,t)&\displaystyle\ge c\varepsilon\int_{-b}^{-a}(t-s)^{\sigma-1}\int_{|z|<\delta}K_{t-s}(x-z)\,dzds \\ [3mm]
&\displaystyle\ge c\varepsilon\int_{t+a}^{t+b} s^{\sigma-1-\frac N2}e^{-\frac{|x|^2}{2s}}\int_{|z|<\delta}e^{-\frac{|z|^2}{2s}}\,dzds\\ [3mm]
&\displaystyle\ge c\varepsilon\delta^N\int_{t+a}^{t+b} s^{\sigma-1-\frac N2}e^{-\frac{|x|^2}{2s}}\,ds\\ [3mm]
&\displaystyle\ge c\varepsilon\delta^Ne^{-\delta^2}(b-a)t^{\sigma-1}K_{t}(x),
\end{array}
$$
provided $t\gg1$.
\end{proof}

\begin{teo}\label{teo-fujita2}
If $p=1+\frac{2\sigma}{N+2(1-\sigma)}$ then for any given data $f$ the solution to problem \eqref{eq.principal} blows up in finite time.
\end{teo}

\begin{proof}
  As we can deduce from the proof of the subcritical case, a sufficient condition to have blow-up if $p=p_*$ for a given data $f$ is
  $$
\int_\Lambda  g(X)\,d\mu(X)\gg 1,
$$
where $g$ is the $\sigma$--caloric extension of $f$. By the semigroup property of problem~\eqref{eq.principal2}, using $U(X,t)$ as new initial datum for any given $t>0$, a sufficient condition is then
$$
\int_\Lambda  U(X,t)\,d\mu(X)\gg 1,\qquad\text{for $t$ large},
$$
see the proof in \cite{Weissler} for the local case.

Integrating the equation in \eqref{eq.principal2} we obtain
$$
\int_\Lambda  U(X,t)\,d\mu(X)=\int_\Lambda  g(X)\,d\mu(X)+\int_0^t\kappa_\sigma^{-1}\int_{\mathbb{R}^N}u^p(x,s)\,dxds.
$$
Introducing now estimate \eqref{eq.homogeneous2} inside the second term, we get
$$
\begin{array}{rl}
\displaystyle\int_\Lambda  U(X,t)\,d\mu(X)
&\displaystyle \ge c\int_{t_0}^t\int_{\mathbb{R}^N}\left(s^{\sigma-1}K_{s}(x)\right)^p\,dxds \\ [3mm]
&\displaystyle \ge c\int_{t_0}^t s^{(\sigma-1-\frac N2)p}\int_{\mathbb{R}^N}e^{-\frac{p|x|^2}{4s}}\,dxds\\ [3mm]
&\displaystyle = c\int_{t_0}^t s^{-1}\,ds\to\infty
\end{array}
$$
for the precise value of $p$.
\end{proof}

We now show that above the critical exponent $p_*$ there exist small global solutions. The existence of such solutions is proved by comparison in the original master equation in \eqref{eq.principal} and the pointwise expression  \eqref{kernel-def} of the nonlocal operator.

\begin{teo}\label{teo-fujita3}
If $p>1+\frac{2\sigma}{N+2(1-\sigma)}$ then the solution to problem \eqref{eq.principal} is global in time provided the datum $f$ is small enough. Moreover the solution tends to zero as $t\to\infty$.
\end{teo}

\begin{proof} We try a supersolution of the form $v(x,t)=\psi(x,t+t_0)$, $\psi(x,t)=t^\nu G(x,t)$ for some $t_0,\nu>0$, where $G$ is given in \eqref{fundamental}.

We write $\psi(x,t)=t^\eta K_t(x)\mathds{1}_{\{t>0\}}$, $\eta>\sigma-1$, and calculate
$$
\begin{array}{l}
\mathcal{M} \psi(x,t)= \\ [3mm]
\quad\displaystyle=\int_{-\infty}^t\int_{\mathbb{R}^N}\left(t^\eta K_t(x)\mathds{1}_{\{t>0\}}-s^\eta K_s(z)\mathds{1}_{\{s>0\}}\right)K_{t-s}(x-z)L(t-s)\,dzds \\ [3mm]
\displaystyle\quad=\frac1{|\Gamma(-\sigma)|}K_t(x)
\int_{-\infty}^t\left(t^{\eta}\mathds{1}_{\{t>0\}}-s^{\eta}\mathds{1}_{\{s>0\}}\right)(t-s)^{-\sigma-1}\,ds\\ [3mm]
\displaystyle\quad=\frac{\Gamma(\eta+1)}{\Gamma(\eta+1-\sigma)}t^{\eta-\sigma}K_t(x)\mathds{1}_{\{t>0\}}.
\end{array}
$$
We have used the semigroup property $K_s*K_{t-s}=K_t$, 
plus the Marchaud derivative of a power. The condition $\eta>\sigma-1$ is critical.

Then $\psi$ is a supersolution of the equation if
$$
ct^{\eta-\sigma-\frac N2}e^{-\frac{|x|^2}{4t}}\ge t^{p(\eta-\frac N2)}e^{-\frac{p|x|^2}{4t}},
$$
for $t>t_0$ large, that is, if
$$
\sigma-1<\eta<\frac N2-\frac\sigma{p-1},\qquad \text{i.e.}\quad p>1+\frac{2\sigma}{N+2(1-\sigma)}.
$$
Then every solution with initial datum $f$ below $\psi$ for $t\le0$ is global in time.

Finally, since $\eta<\frac N2$, the solution decreases in time to zero.
\end{proof}

\

We can obtain a sharper result in this supercritical Fujita case provided we are in the subcritical (fractional) Sobolev range.

\begin{teo}\label{teo-abovefujita}
If $1+\frac{2\sigma}{N+2(1-\sigma)}<p<\frac{N+2\sigma}{N-2\sigma}$, $N>2\sigma$, then any global  solution to problem \eqref{eq.principal}  tends to zero as $t\to\infty$.
\end{teo}

The proof is based on the so-called concavity argument of Levine in terms of the energy defined in \eqref{energy}, see \cite{Levine} for the local case, together with some nonexistence of stationary solutions result. Here it can be checked with trivial computations, that if the energy is initially negative, then there exists appropriate positive constants $\alpha$ and $A$ such that  the function
$$
H(t)=\int_0^t\int_\Lambda U^2(X,t)\,d\mu(X)+A
$$
satisfies
$$
H,\;H'\ge0,\qquad (H^{-\alpha})''\le0.
$$
This implies that $H$ cannot be global in time. We then have the following result.
\begin{teo}\label{Levine}
If $I_U(t)<0$ for some $t\ge0$  then the solution to problem \eqref{eq.principal2} blows up in finite time.
\end{teo}

\begin{proof}[Proof of Theorem \ref{teo-abovefujita}]
If $U$ is a global solution then $I_U(t)\ge0$ for every $t\ge0$. Since also $I_U(t)\le I_U(0)$ we get, through a standard Lyapunov argument, that $\lim\limits_{t\to\infty}U(\cdot,t)=Z$,  a stationary solution to the equation in~\eqref{eq.principal2}. If $p<\frac{N+2\sigma}{N-2\sigma}$, $N>2\sigma$ we conclude that $Z\equiv0$, see \cite{BCdPS-13,CLO-05}.
\end{proof}

\

\section{Blow-up rate}\label{sec-rate}
In this Section we prove the blow-up rate below the Fujita exponent, using a rescaling technique inspired in the work \cite{FilaSouplet01}, see also \cite{Hu90,FdP21}, and thanks to the fact that all solutions blow up in that range of parameters.

\begin{proof}[Proof of Theorem \ref{teo-rate}] The lower blow-up rate in \eqref{rate}
is clear from the existence of the explicit solution
$$
z(t)=c(T-t)^{-\frac\sigma{p-1}},\qquad c^{p-1}=\frac{\Gamma(\sigma p/(p-1)}{\Gamma(\sigma/(p-1)},
$$
plus the strong maximum principle (two solutions having the same blow-up time must intersect). Let us then prove the upper estimate.

Define
$$
m(t)=\max_{\mathbb{R}^N\times (-\infty,t]} u(x,\tau)
$$
and consider an increasing sequence of times
$$
t_{j+1}=\sup\{t\in(t_j,T)\,:\, m(t)=2m(t_j)\}.
$$
Observe that $m(t)$ is a nondecreasing function and in general $\|u(\cdot, t)\|_\infty\le m(t)$
($u$ is not necessarily increasing in time).
However, for this sequence we have
$\|u(\cdot, t_j)\|_\infty\le m(t_j)$. Therefore we can take $x_j\in\mathbb{R}^N$ such that
$$
u(x_j,t_j)= m(t_j).
$$

Assume
\begin{equation}\label{eq.times}
\lim_{j\to\infty}(t_{j+1}-t_j){m^{\nu}(t_j)}=\infty,\qquad \nu=\frac{p-1}\sigma,
\end{equation}
since on the contrary we would easily obtain the desired rate. In fact, if \eqref{eq.times} does
not hold,
$$
t_{j+1} - t_j \le cm^{-\nu}(t_j) = c2^{j(1-p)}m^{-\nu}(t_0),
$$
so performing the sum we get
$$
T - t_0 \le cm^{-\nu}(t_0)\sum_{j=0}^\infty
2^{-j\nu} = cm^{-\nu}(t_0) \le c\|u(\cdot, t_0)\|_\infty^{-\nu}
$$
that is, \eqref{rate}.

Define now the functions
$$
\varphi_j(y,s)=\frac1{m_j} u(m_j^{-\nu }y+ x_j,m_j^{-2\nu}s+ t_j),
$$
for
$$
y\in \mathbb R^N, \qquad s\in I_j= (- \infty, (T- t_j)m_j^{2\nu}),
$$
where $m_j=m(t_j)$. Notice that  $\varphi_j$ is a solution of the equation
$$
\mathcal{M}(\varphi_j)=\varphi_j^p,\qquad  (y,s)\in \mathbb R^N\times I_j,
$$
which satisfies
$$
\begin{array}{l}
\varphi_j(0,0)=1,\\
\varphi_j(y,s)\le 2 \quad \mbox{ in } \mathbb R^N\times (- \infty,(t_{j+1}-t_j)m_j^{2\nu}).
\end{array}
$$
Moreover, thanks to \eqref{eq.times} we have that
$$
I_j\to \mathbb R\qquad \mbox{as } j\to\infty.
$$
Since the sequence $\varphi_j$ is uniformly bounded, we know that $\varphi_j\in C^{\varepsilon}_p$ locally for some $\varepsilon>0$. Also, since $\varphi_j(0,0)=1$ we have a uniform nontrivial lower bound for every $\varphi_j$. That is,
$$
\varphi_j(y,s)\ge g(y,s)\ge0,
$$
for instance $g(y,s)=[1-c(|y|^2+|s|)^{\varepsilon/2}]_+$.

Let $h(x,t)$ be the solution of
$$
\left\{
\begin{array}{ll}
\mathcal{M} h=h^p,\qquad & (y,s)\in \mathbb{R}^N \times(0,(t_{j+1}-t_j)m_j^{p-1}),\\
h(y,s)=g(y,s),\qquad & (y,s)\in \mathbb{R}^N \times(-\infty,0).
\end{array}\right.
$$
We end by applying Theorem \ref{teo-BUP} so that $h$ blows up in finite time $S$ independent of $j$, and so does each $\varphi_j$. Take $j$  large so that $(t_{j+1}-t_j)m_j^{2\nu}>S$ to get a contradiction with the fact that $\varphi_j$ is bounded.
\end{proof}

\section*{Appendix}

We prove here the existence of a unique solution to our problem \eqref{eq.principal}. We introduce the following notation: the parabolic H\"older space $C_p^\gamma(\Omega)$, $\Omega\subset\mathbb{R}^{N+1}$, $0<\gamma<1$, is defined by
\begin{equation}\label{space-C-gamma}
v\in C_p^\gamma(\Omega)\;\Leftrightarrow\; v\in L^\infty(\Omega),\;|v(x,t)-v(z,s)|\le C\left(|x-z|^2+|t-s|\right)^{\gamma/2},
\end{equation}
for every $(x,t),(z,s)\in\Omega$. We also consider the norm
$$
\|v\|_{C_p^\gamma}=\|v\|_\infty+\inf\{C>0\;:\;\text{estimate } \eqref{space-C-gamma} \text{ holds}\}.
$$
The analogous spaces for $\gamma\ge1$ can also be considered.

The operator $\mathcal{M}$ is well defined for functions in this space when $\gamma>2\sigma$, thus compensating the singularity of the kernel at the space-time origin.

\begin{proof}[Proof of Theorem \ref{teo-existence}]

We first consider the linear problem
\begin{equation}\label{eq.general2}
\mathcal{M} u=F,\qquad (x,t)\in Q_T=\mathbb{R}^N\times(-\infty,T).
\end{equation}
This problem is studied in \cite{StingaTorrea17}.
If $F\in C_p^\varepsilon(Q_T)$, for some $\varepsilon>0$, and decays at $t=-\infty$ faster than $|t|^{-\sigma}$, then the solution $u\in C_p^{\varepsilon+2\sigma}(Q_T)$ to problem \eqref{eq.general2} is given by the formula
\begin{equation*}\label{repres-linear}
u(x,t)=\int_{-\infty}^t\int_{\mathbb{R}^N}F(z,s)G(x-z,t-s)\,dzds,
\end{equation*}
where $G$ is the Green function \eqref{fundamental}. The proof performed in \cite{StingaTorrea17} uses Fourier transform, and then $F\in L^2(\mathbb{R}^{N+1})$ is moreover required, but it holds also under our less restrictive conditions.

As a second step we consider the linear problem with memory
\begin{equation}\label{eq.general3}
\left\{
\begin{array}{ll}
\mathcal{M} u=h,\qquad & \text{for } (x,t)\in Q^+_T=\mathbb{R}^N\times(0,T),\\
u=f,\qquad & \text{for } (x,t)\in Q_0,
\end{array}\right.
\end{equation}
In order to apply the previous result we take
\begin{equation*}\label{representation-F2}
F=\begin{cases}
  h&\text{if } t>0, \\
  \mathcal{M}f&\text{if } t\le0.
\end{cases}
\end{equation*}
We need $h\in C_p^{\varepsilon}(Q_T^+)$ and $f\in C_p^{\varepsilon+2\sigma}(Q_0)$ satisfying the decay condition $|\mathcal{M}f(t)|\le c|t|^{-\sigma-\varepsilon}$ for $t\le-1$.

Finally a standard fixed point argument allows to consider the nonlinear problem \eqref{eq.principal}.
Let $T_0>0$, $\varepsilon>0$, $K>0$ be given and let $f\in C_p^{\sigma+\varepsilon}(Q_0)$ satisfying condition \eqref{H}. We consider the space
$$
E_{T_0}=\{v\in C_p^{\sigma+\varepsilon}(Q_{T_0}) \;:\; v\ge0,\;\|v\|_{C_p^{\sigma+\varepsilon}}\le K\|f\|_{C_p^{\sigma+\varepsilon}}\}.
$$
Define also the operator $R$ on $E_{T_0}$
\begin{equation}\label{operator}
  Rv=R_-f+R_+v^p,
\end{equation}
where
$$
\begin{array}{l}
\displaystyle R_-f(x,t)=\int_{-\infty}^0 \int_{\mathbb{R}^N}\mathcal{M}f(z,s)G(x-z,t-s)\,dzds, \\ [3mm]
\displaystyle R_+v(x,t)=\int_{0}^t\int_{\mathbb{R}^N}v(z,s)G(x-z,t-s)\,dzds.
\end{array}
$$
The function $Rv$ solves problem \eqref{eq.general3} with $h=v^p$.
The estimates obtained in \cite{StingaTorrea17} for problem \eqref{eq.general2}, or which is the same for problem \eqref{eq.general3}, allow to check that this operator is well defined in $E_{T_0}$ if $K>0$ is large and contractive if $T_0$ is small enough. Therefore there exists a fixed point, which is a mild solution to our problem. Uniqueness and existence for a maximal time interval $T>0$ are also easy. Regularity follows by bootstrapping the results of \cite{StingaTorrea17} starting from a bounded reaction $u^p$ for $0<t<T-\varepsilon$, thus getting $2\sigma$ more regularity in each step, so the solution possesses the same regularity as the memory data $f$.
\end{proof}

\section*{Acknowledgments}

Work supported by the Spanish project  PID2020-116949GB-I00. The first author was also supported by Grupo de Investigaci\'on UCM 920894.

\end{document}